\magnification\magstep1
\input  amssym.tex

\overfullrule0pt

\def\sqr#1#2{{\vcenter{\hrule height.#2pt              %qed
     \hbox{\vrule width.#2pt height#1pt\kern#1pt
     \vrule width.#2pt}
     \hrule height.#2pt}}}
\def\square{\mathchoice\sqr{5.5}4\sqr{5.0}4\sqr{4.8}3\sqr{4.8}3}
\def\qed{\hskip4pt plus1fill\ $\square$\par\medbreak}

\centerline{\bf No smooth Julia sets for polynomial}

\centerline{\bf  diffeomorphisms of ${\Bbb C}^2$ with positive entropy}

\bigskip \centerline{Eric Bedford and Kyounghee Kim}

\bigskip
\noindent{\bf \S0.  Introduction.}  There are several reasons why the polynomial diffeomorphisms of ${\Bbb C}^2$ form an interesting family of dynamical systems. One of these is the fact that there are connections with two other areas of dynamics: polynomial maps of ${\Bbb C}$ and diffeomorphisms of ${\Bbb R}^2$, which have each received a great deal of attention.  The question arises whether, among the polynomial diffeomorphisms of ${\Bbb C}^2$, are there maps with the special status of having smooth Julia sets?  Here we show that is not the case.

More generally, we consider a holomorphic mapping $f:X\to X$ of a complex manifold $X$.  The Fatou set of $f$ is defined as the set of points $x\in X$ where the iterates $f^n:=f\circ\cdots\circ f$ are locally equicontinuous.  If $X$ is not compact, then in the definition of equicontinuity, we consider the one point compactification of $X$; in this case, a sequence which diverges uniformly to infinity is equicontinuous.   By the nature of equicontinuity, the dynamics of $f$ is regular on the Fatou set.  The Julia set is defined as the complement of the Fatou set, and this is where any chaotic dynamics of $f$ will take place.  The first nontrivial case is where $X={\Bbb P}^1$ is the Riemann sphere, and in this case Fatou (see [M1]) showed that if the Julia set $J$ is a smooth curve, then either $J$ is the unit circle, or $J$ is a real interval.  If $J$ is the circle, then $f$ is equivalent to $z\mapsto z^d$, where $d$ is an integer with $|d|\ge2$;  if $J$ is the interval, then $f$ is equivalent to a Chebyshev polynomial.  These maps with smooth $J$ play special roles, and this sparked our interest to look for smooth Julia sets in other cases.

Here we address the case where $X={\Bbb C}^2$, and $f$ is a polynomial automorphism, which means that $f$ is biholomorphic, and the coordinates are polynomials.  Since $f$ is invertible, there are two Julia sets: $J^+$ for iterates in forward time, and $J^-$ for iterates in backward time.   Polynomial automorphisms have been classified by Friedland and Milnor [FM]; every such automorphism is conjugate to a map which is either affine or elementary, or it belongs to the family ${\cal H}$.   The affine and elementary maps have simple dynamics, and $J^\pm$ are (possibly empty) algebraic sets (see [FM]).

Thus we will restrict our attention to the maps in ${\cal H}$, which are finite compositions $f:=f_k\circ\cdots\circ f_1$, where each $f_j$ is a generalized H\'enon map, which by definition has the form $f_j(x,y)=(y, p_j(y)-\delta_j x)$, where $\delta_j\in{\Bbb C}$ is nonzero, and $p_j(y)$ is a monic polynomial of degree $d_j\ge2$.  The degree of $f$ is $d:=d_1\cdots d_k$, and the complex Jacobian of $f$ is $\delta:=\delta_1\cdots\delta_k$.    In [FM] and [Sm] it is shown that the topological entropy of $f$ is $\log d>0$.  The dynamics of such maps is complicated and has received much study, starting with the papers [H], [HO1], [BS1] and [FS].   

For maps in ${\cal H}$, we can ask whether $J^+$ can be a manifold.  For any saddle point $q$, the stable manifold $W^s(q)$ is a Riemann surface contained in  $J^+$.  Thus $J^+$ would have to have real dimension at least two.  However, $J^+$ is also the support of a positive, closed current $\mu^+$ with continuous potential, and such potentials cannot be supported on a Riemann surface (see [BS1, FS]).   On the other hand, since $J^+=\partial K^+$ is a boundary, it cannot have interior.  Thus dimension 3 is the only possibility for $J^+$ to be a manifold.  In fact, there are examples of $f$ for which $J^+$ has been shown to be a topological 3-manifold (see [FS], [HO2],  [Bo], [RT]).   
%On the other hand, Forn\ae ss and Sibony [FS] have noted that $J^+$ cannot be smooth for a generic element of ${\cal H}$.

%
%We say that $X\subset{\Bbb C}^2$ is a {\it 3 dimensional manifold-with-boundary} to mean that $X$ is closed, and each point $q\in X$ has a neighborhood in $X$ which is homeomorphic to either ${\Bbb R}^3$ or $\{(x_0,x_1,x_2)\in{\Bbb R}^3:x_0\ge0\}$.  The boundary of $X$ is the 2 real dimensional set given locally as $\{x_0=0\}$.  
The purpose of this paper is to prove the following:
 \proclaim   Theorem.  For any polynomial automorphism of ${\Bbb C}^2$ of positive entropy, neither $J^+$ nor $J^-$ is smooth of class $C^1$, in the sense of manifold-with-boundary.
  
\noindent    We may interchange the roles of $J^+$ and $J^-$ by replacing $f$ by $f^{-1}$, so there is no loss of generality if we consider only $J^+$.

In an Appendix, we discuss the non-smoothness of the related sets $J$, $J^*$, and $K$.
\medskip
\noindent{\it Acknowledgment.}  We wish to thank Yutaka Ishii and Paolo Aluffi for helpful conversations on this material.

\medskip\noindent{\bf \S1. No boundary. }  Let us start by showing that if $J^+$ is a $C^1$ manifold-with-boundary, then the boundary is empty.  Recall that if $J^+$ is $C^1$, then for each $q_0\in J^+$ there is a neighborhood $U\ni q_0$ and $r,\rho\in C^1(U)$ with $dr\wedge d\rho\ne0$ on $U$,  such that $U\cap J^+=\{r=0, \rho\le0\}$.  If $J^+$ has boundary, it is given locally by $\{r=\rho=0\}$.   For $q\in J^+$, the tangent space $T_qJ^+$ consists of the vectors that annihilate $dr$.  This contains the subspace $H_q\subset T_qJ^+$, consisting of the vectors that annihilate $\partial r$.  $H_q$ is the unique complex subspace inside $T_qJ^+$, so if $M\subset J^+$ is a complex submanifold, then $T_qM=H_q$.

We start by showing that if $J^+$ is  $C^1$, then it carries a Riemann surface lamination.
\proclaim Lemma 1.1.  If $J^+$ is $C^1$ smooth, then $J^+$ carries a Riemann surface foliation ${\cal R}$ with the property that if $W^s(q)$ is the stable manifold of a saddle point $q$, then $W^s(q)$ is a leaf of ${\cal R}$.  If $J^+$ is a $C^1$ smooth manifold-with-boundary, then ${\cal R}$ extends to a Riemann surface lamination of $J^+$.  In particular, any boundary component is a leaf of ${\cal R}$.

\noindent{\it Proof. }   Given $q_0\in J^+$, let us choose holomorphic coordinates $(z,w)$ such that  $dr(q_0) = dw$.  We work in a small neighborhood which is a bidisk $\Delta_\eta\times\Delta_\eta$.  We may choose $\eta$ small enough that $|r_z/r_w|<1$.  In the $(z,w)$-coordinates, the tangent space $H_q$ has slope less than 1 at every point $\{|z|,|w|<\eta\}$.   Now let $\hat q$ be a saddle point, and let $W^s(\hat q)$ be the stable manifold, which is a complex submanifold of ${\Bbb C}^2$, contained in $J^+$.  Let $M$ denote a connected component of $W^s(\hat q)\cap(\Delta_\eta\times\Delta_{\eta/2})$.   Since the slope is $<1$, it follows that there is an analytic function $\varphi:\Delta_\eta\to\Delta_\eta$ such that $M\subset\Gamma_\varphi:=\{(z,\varphi(z)):z\in \Delta_\eta\}$.  Let $\Phi$ denote the set of all such functions $\varphi$.  Since a stable manifold can have no self-intersections, it follows that if $\varphi_1,\varphi_2\in\Phi$, then either $\Gamma_{\varphi_1}=\Gamma_{\varphi_2}$ or $\Gamma_{\varphi_1}\cap\Gamma_{\varphi_2}=\emptyset$.   Now let $\hat\Phi$ denote the set of all normal limits (uniform on compact subsets of $\Delta_\eta$) of elements of $\Phi$.  We note that by Hurwitz's Theorem, the graphs $\Gamma_\varphi$, $\varphi\in\hat\Phi$ have the same pairwise disjointness property.  Finally, by [BS2], $W^s(q_0)$ is dense in $J^+$, so the graphs $\Gamma_\varphi$, $\varphi\in\hat\Phi$ give the local Riemann surface lamination.

If $q_1$ is another saddle point, we may follow the same procedure and obtain a Riemann surface lamination whose graphs are given locally by $\varphi\in\hat\Phi_1$.  However, we have seen that the tangent space to the foliation at a point $q$ is given by $H_q$.  Since these two foliations have the same tangent spaces everywhere, they must coincide.

We have seen that all the graphs are contained in $J^+$, so if $J^+$ has boundary, then the boundary must coincide locally with one of the graphs.
\qed

We will use the observation that $K^+\subset\{(x,y)\in{\Bbb C}^2:|y|>\max(|x|,R)\}$.  Further, we will use the Green function $G^-$ which has many properties, including:  
\item{$(i)$} $G^-$ is pluri-harmonic on $\{G^->0\}$, 
\item{$(ii)$} $\{G^-=0\}=K^-$, 
\item{$(iii)$} $G^-\circ f = d^{-1}G^-$.  

\noindent Further, the restriction of $G^-$ to $\{|y|\le \max(|x|,R)\}$ is a proper exhaustion.
\proclaim Lemma 1.2.  Suppose that $J^+$ is a $C^1$ smooth manifold-with-boundary, and $M$ is a component of the boundary of $J^+$.  Then $M$ is a closed Riemann surface, and $M\cap K\ne\emptyset$.

\noindent{\it Proof.  }  We consider the restriction $g:=G^-|_{M}$.  If $M\cap K=\emptyset$, then $g$ is harmonic on $M$.  On the other hand, $g$ is a proper exhaustion of $M$, which means that $g(z)\to\infty$ as $z\in M$ leaves every compact subset of $M$.  This means that $g$ must assume a minimum value at some point of $M$, which would violate the minimum principle for harmonic functions.
\qed

\proclaim Lemma 1.3.    Suppose that $J^+$ is a $C^1$ smooth manifold-with-boundary, then the boundary is empty.

\noindent{\it Proof.  }  Let $M$ be a component of the boundary of $J^+$.  By Lemma  1.2, $M$ must intersect $\Delta_R^2$.  Since $J^+$ is $C^1$, there can only finitely many boundary components of $J^+\cap \Delta_R^2$.  Thus there can be only finitely many components $M$, which must be permuted by $f$.  If we take a sufficiently high iterate $f^N$, we may assume that $M$ is invariant.  Now let $h:=f^N|_{M}$ denote the restriction to $M$.  We see that $h$ is an automorphism of the Riemann surface $M$, and the iterates of all points of $M$ approach $K\cap M$ in forward time.  It follows that $M$ must have a fixed point $q\in M$, and $|h'(q)|<1$.  
The other multiplier of $Df$ at $q$ is $\delta/h'(q)$.  

We consider three cases.  First, if $|\delta/h'(q)|>1$, then $q$ is a saddle point, and $M=W^s(q)$.  On the other hand, by [BS2], the stable manifold of a saddle points is dense in $J^+$, which makes it impossible for $M$ to be the boundary of $J^+$.  This contradiction means that there can be no boundary component $M$.  

The second case is $|\delta/h'(q)|<1$.  This case cannot occur because the multipliers are less than 1, so $q$ is a sink, which means that $q$ is contained in the interior of $K^+$ and not in $J^+$.

The last case is where  $|\delta/h'(q)|=1$.  In this case, we know that $f$ preserves $J^+$, so $Df$ must preserve $T_q(J^+)$.  This means that the outward normal to $M$ inside $J^+$ is preserved, and thus the second multiplier must be +1.  It follows that $q$ is a semi-parabolic/semi-attracting fixed point.  It follows that $J^+$ must have a cusp at $q$ and cannot be $C^1$ (see Ueda [U] and Hakim [Ha]).
\qed

\medskip\noindent{\bf \S2.  Maps that do not decrease volume.}  We note the following topological result (see Samelson [S] for an elegant proof): {\sl If $M$ is a smooth 3-manifold (without boundary) of class $C^1$ in ${\Bbb R}^4$, then it is orientable.}   This gives:

\proclaim Proposition 2.1.  For any $q\in M$, there is a neighborhood $U$ about $q$ so that $U-M$ consists of two components ${\cal O}_1$ and ${\cal O}_2$, which belong to different components of ${\Bbb R}^4-M$.  

\noindent{\it Proof. }  Suppose that ${\cal O}_1$ and ${\cal O}_2$ belong to the same component of ${\Bbb R}^4-M$.  Then we can construct a simple closed curve $\gamma\subset{\Bbb R}^4$ which crosses $M$ transversally at $q$ and has no other intersection with $M$.  It follows that the (oriented) intersection  is $\gamma\cdot M=1$ (modulo 2).  But the oriented intersection modulo 2 is a homotopy invariant (see [M2]), and $\gamma$ is contractible in ${\Bbb R}^4$, so we must have $\gamma\cdot M=0$ (modulo 2).  \qed

\proclaim Corollary 2.2.  If $J^+$ is $C^1$ smooth, then $f$ is an orientation preserving map of $J^+$.

\noindent{\it Proof. }  $U^+:={\Bbb C}^2-K^+$ is a connected (see [HO1]) and thus it is a component of ${\Bbb C}^2-J^+$.  Since $f$ preserves $U^+$, it also preserves the orientation of $J^+$, which is $\pm \partial U^+$.  \qed

We recall the following result of Friedland and Milnor:
\proclaim Theorem [FM].  If $|\delta|>1$, then $K^+$ has zero Lebesgue volume, and thus $J^+=K^+$.  If $|\delta|=1$, then ${\rm int}(K^+)={\rm int}(K^-) = {\rm int}(K)$.  In particular, there exists $R$ such that $J^+=K^+$ outside $\Delta_R^2$.

\noindent{\it Proof of  Theorem in the case $|\delta|\ge1$. }   Let $q\in J^+$ be a point outside $\Delta^2_R$, as in the Theorem above.  Then near $q$ there must be a component ${\cal O}$, which is distinct from $U^+={\Bbb C}^2-K^+$.   Thus ${\cal O}$ must belong to the interior of $K^+$.  But by the Theorem above, the interior of $K^+$ is not near $q$.  \qed

\medskip\noindent{\bf \S3.  Volume decreasing maps. }    Throughout this section, we continue to suppose that $J^+$ is $C^1$ smooth, and in addition we suppose that $|\delta|<1$.   For a point $q\in J^+$, we let $T_q:=T_q(J^+)$ denote the real tangent space to $J^+$.  We let $H_q:=T_q\cap i T_q$ denote the unique (one-dimensional) complex subspace inside $T_q$.  Since $J^+$ is invariant under $f$, so is $H_q$, and we let $\alpha_q$ denote the multiplier of $D_qf|_{H_q}$.  

\proclaim Lemma 3.1.  Let $q\in J^+$ be a fixed point.  There is a $D_qf$-invariant subspace $E_q\subset T_q({\Bbb C}^2)$ such that $H_q$ and $E_q$ generate $T_q$.   We denote the multiplier of $D_qf|_{E_q}$ by $\beta_q$.  Thus $D_qf$ is linearly conjugate to the diagonal matrix with diagonal elements $\alpha_q$ and $\beta_q$.  Further, $\beta_q\in{\Bbb R}$, and $\beta_q>0$.

\noindent{\it Proof. }  We have identified an eigenvalue $\alpha_q$ of $D_qf$.   If $D_qf$ is not diagonalizable, then it must have a Jordan canonical form $\pmatrix{\alpha_q & 1\cr 0& \alpha_q}$.    The determinant is $\alpha_q^2=\delta$, which has modulus less than 1.  Thus $|\alpha_q|<1$, which means that $q$ is an attracting fixed point and thus in the interior of $K^+$, not in $J^+$.  Thus $D_qf$ must be diagonalizable, which means that $H_q$ has a complementary invariant subspace $E_q$.  
Since $E_q$ and $T_q$ are invariant under $D_qf$, the real subspace $E_q\cap T_q\subset E_q$ is invariant, too.  Thus $\beta_q\in{\Bbb R}$.  By Corollary~2.2,  $D_qf$ will preserve the orientation of $T_q$, and so $\beta_q>0$.  \qed

Let us recall the Riemann surface foliation of $J^+$ which was obtained in Lemma 1.1.  For $q\in J^+$, we let $R_q$ denote the leaf of ${\cal R}$ containing $q$.  If $q$ is a fixed point, then $f$ defines an automorphism $g:=f|_{R_q}$ of the Riemann surface $R_q$.  Since $R_q\subset K^+$, we know that the iterates of $g^n$ are bounded in a complex disk $q\in \Delta_q\subset R_q$.  Thus the derivatives $(Dg)^n= D(g^n)$ are bounded at $q$.  We conclude that $|\alpha_q|=|D_q(g)|\le 1$.  If $|\alpha_q|=1$, then $\alpha_q$ is not a root of unity.  Otherwise $g$ is an automorphism of $R_q$ fixing $q$, and $Dg^n(q) =1$ for some $n$.  It follows that $g^n$ must be the identity on $R_q$.  This means that $Rq$ would be a curve of fixed points for $f^n$, but by [FM] all periodic points of $f$ are isolated, so this cannot happen.
\proclaim Lemma 3.3.  If $q\in J^+$ is a fixed point, then $q$ is a saddle point, and  $\alpha_q=\delta/d$, and $\beta_q=d$.

\noindent {\it Proof. }  First we claim that $|\alpha_p|<1$.  Otherwise, we have $|\alpha_q|=1$, and by the discussion above, this means that $\alpha_q$ is not a root of unity.  Thus the restriction $g = f|_{R_q}$ is an irrational rotation.  Let $\Delta\subset R_q$ denote a $g$-invariant disk containing $q$.  Since $|\delta|=|\alpha_q\beta_q|=|\beta_q|$ has modulus less than 1, we conclude that $f$ is normally attracting to $\Delta$, and thus $q$ must be in the interior of $K^+$, which contradicts the assumption that $q\in J^+$.

Now we have $|\alpha_q|<1$, so if $|\beta_q|=1$, we have $\beta_q=1$, since $\beta_q$ is real and positive.  This means that $q$ is a semi-parabolic, semi-attracting fixed point for $f$.  We conclude by Ueda [U] and Hakim [Ha] that $J^+$ has a cusp at $q$ and thus is not smooth.  Thus we conclude that $|\beta_q|>1$, which means that $q$ is a saddle point.

Now since $E_q$ is transverse to $H_q$, it follows that $W^u(q)$ intersects $J^+$ transversally, and thus $J^+\cap W^u(q)$ is $C^1$ smooth.  Let us consider the uniformization
$$\phi:{\Bbb C}\to W^u(q)\subset{\Bbb C}^2, \ \ \ \phi(0)=q, \ \ \ f\circ\phi(\zeta)=\phi(\lambda^u\zeta) $$
The pre-image $\tau:=\phi^{-1}(W^u(q)\cap J^+)\subset{\Bbb C}$ is a $C^1$ curve passing through the origin and invariant under $\zeta\mapsto\lambda^u\zeta$.   It follows that $\lambda^u\in {\Bbb R}$, and $\tau$ is a straight line containing the origin.  Further, $g^+:=G^+\circ\phi$ is harmonic on ${\Bbb C}-\tau$, vanishing on $\tau$, and satisfying $g^+(\lambda^u\zeta)=d\cdot g^+(\zeta)$.  Since $\tau$ is a line, it follows that $g^+$ is piecewise linear, so we must have $\lambda^u=\pm d$.  Finally, since $f$ preserves orientation, we have $\lambda^u=d$.
\qed

\proclaim Lemma 3.4.  There can be at most one fixed point in the interior of $K^+$.  There are at least $d-1$  fixed points are contained in $J^+$, and at each of these fixed points, the differential $Df$ has multiplier of $d$.

\noindent{\it Proof. }  Suppose that $q$ is a fixed point in the interior of $K^+$.  Then $q$ is contained in a recurrent Fatou domain $\Omega$, and by [BS2], $\partial\Omega=J^+$.  If there is more than one fixed point in the interior of $K^+$, we would have $J^+$ simultaneously being the boundary of more than one domain, in addition to being the boundary of $U^+={\Bbb C}^2-K^+$.  This is not possible if $J^+$ is a topological submanifold of ${\Bbb C}^2$.  

By [FM] there are exactly $d$ fixed points, counted with multiplicity.  By Lemma~3.3, the fixed points in $J^+$ are of saddle type, so they have multiplicity 1.  Thus there are at least $d-1$ of them.  \qed

\medskip\noindent{\bf \S4.  Fixed points with given multipliers.}  If $q=(x,y)$ is a fixed point for $f=f_n\circ\cdots\circ f_1$, then we may represent it as a finite sequence $(x_j,y_j)$ with $j\in{\Bbb Z}/n{\Bbb Z}$, subject to the conditions $(x,y)=(x_1,y_1)=(x_{n+1},y_{n+1})$ and $f_j(x_j,y_j)=(x_{j+1},y_{j+1})$.  Given the form of $f_j$, we have $x_{j+1}=y_j$, so we may drop the $x_j$'s from our notation and write $q=(y_n,y_1)$.  We identify this point with the sequence $\hat q=(y_1,\dots,y_n)\in{\Bbb C}^n$, and we define the polynomials
$$\eqalign{\varphi_1&:= p_1(y_1)-\delta_1 y_n - y_2  \cr
\varphi_2&:=p_2(y_2) - \delta_2 y_1 - y_3\cr
&\dots   \dots\dots   \cr
\varphi_n &:=  p_n(y_n) - \delta_n y_{n-1} - y_1}$$
The condition to be a fixed point is that $\hat q=(y_1,\dots,y_n)$  belongs to the zero locus $Z(\varphi_1,\dots,\varphi_n)$ of the $\varphi_i$'s.
We define $q_i(y_i):=p_i(y_i) - y_i^{d_i} $  and $Q_i := q_i(y_i) - y_{i+1}-\delta_i y_{i-1}$, so
$$\varphi_i = y_i^{d_i} + q_i(y_i) - y_{i+1} - \delta_i y_{i-1} = y_i^{d_1} + Q_i \eqno(*)$$
Since $p_j$ is monic, the degrees of $q_i$ and $Q_i$ are $\le d_i-1$.

By the Chain Rule, the differential of $f$ at $q=(y_n,y_1)$ is given by
$$Df(q)=\pmatrix{0&1\cr -\delta_n& p'_n(y_n)} \cdots  \pmatrix{0&1\cr -\delta_1& p'_n(y_1)}$$
We will denote this by
$ M_n=M_n(y_1,\dots,y_n):= \pmatrix{m_{11}^{(n)} & m_{12}^{(n)}\cr m_{21}^{(n)} & m_{22}^{(n)} }$.  

We consider special monomials in  $p'_j=p'_j(y_j)$ which have the form $(p')^L:= p'_{\ell_1}\cdots p'_{\ell_s}$, with  $L=\{\ell_1,\dots,\ell_s\}\subset\{1,\dots,n\}$.    Note that the factors $p'_{\ell_i}$ in $(p')^L$ are distinct.  Let us use the notation $|L|$ for the number of elements in $L$, and  $H_{\bf m}$ for the linear span of $\{(p')^L: |L|= m-2k, 0\le k\le  {n/ 2} \}$.
With this notation, {\bf m} indicates the maximum number of factors of $p'_j$ in any monomial, and in every case the number of factors differs from {\bf m} by an even number.
\proclaim Lemma 4.1.  The entries of $M_n$:
\item{(1)} $m_{11}^{(n)}$ and $m_{22}^{(n)}- p'_1(y_1)\cdots p'_n(y_n)$  both belong to $ H_{\bf n-2}$.
\item{(2)}  $m_{12}^{(n)}, m_{21}^{(n)} \in H_{\bf n-1}$.

\noindent{\it Proof. }  We proceed by induction.  The case $n=1$ is clear.  If $n=2$, 
$$M_2=\pmatrix{0&1\cr -\delta_2 & p'_2} \pmatrix{0 & 1\cr -\delta_1 & p'_1} = \pmatrix{-\delta_1 & p'_1\cr -\delta_1p_2'& p'_1p_2'-\delta_2}$$
which satisfies (1) and (2).  For $n>2$, we have
$$M_n = \pmatrix{0&1\cr  -\delta_n & p'_n} M_{n-1} = \pmatrix{m_{21}^{(n-1)} & m_{22}^{(n-1)}\cr  -\delta_n m_{11}^{(n-1)} + m_{21}^{(n-1)} p'_n & -\delta_n m_{12}^{(n-1)} + p'_n m_{22}^{(n-1)}}$$
which gives (1) and (2) for all $n$.  \qed

The condition for $Df$ to have a multiplier $\lambda$ at $q$ is
$\Phi(\hat q)=0$, where
$$\Phi = \det \left(M_n-\pmatrix{\lambda & 0\cr 0& \lambda} \right) $$
\proclaim Lemma 4.2.  $\Phi-p'_1(y_1)\cdots p'_n(y_n) \in H_{\bf n-2}$.

\noindent{\it Proof. }  The formula for the determinant gives
$$\Phi = \lambda^2 - \lambda {\rm Tr}(M_n) + {\rm det}(M_n) = \lambda^2 - \lambda (m_{11}^{(n)} + m_{22}^{(n)}) + \delta$$
since $\delta$ is the Jacobian determinant of $Df$.  The Lemma now follows from Lemma 4.1.
\qed

The degree of the monomial $y^a:= y_1^{a_1}\cdots y_n^{a_n}$ is ${\rm deg}(y^a)=a_1+\cdots+ a_n$.   We will use the {\it graded lexicographical order} on the monomials in $\{y_1,\dots,y_n\}$.  That is, $y^a> y^b$ if either ${\rm deg}(y^a)>{\rm deg}(y^b) $, or if ${\rm deg}(y^a)={\rm deg}(y^b) $ and $a_i>b_i$, where $i = \min\{1\le j\le n: a_j\ne b_j\}$.  If $f\in{\Bbb C}[y_1,\dots,y_n]$, we denote $LT(f)$ for the leading term of $f$, $LC(f)$ for the leading coefficient, and $LM(f)$ for the leading monomial.
\proclaim Lemma 4.3.  With the graded lexicographical order, $G:=\{\varphi_1,\dots,\varphi_n\}$ is a Gr\"obner basis.

\noindent{\it Proof. }   We will use Buchberger's Algorithm (see [CLO, Chapter 2]).  For each $i=1,\dots,n$, $LT(\varphi_i)= LM(\varphi_i)=y_i^{d_i}$, so for $i\ne j$, the least common multiple of the leading terms is $L.C.M. = y_i^{d_i}y_j^{d_j}$.  The $S$-polynomial is
$$S(\varphi_i,\varphi_j) := {L.C.M. \over LM(\varphi_j)}\varphi_i-{L.C.M.\over LM(\varphi_i)}\varphi_j = y_j^{d_j}Q_i-y_i^{d_i}Q_j = \varphi_jQ_i-Q_j\varphi_i$$
where we use the $Q_j$ from (4.1) and cancel terms.  Now let $\mu_i:={\rm deg}(Q_i)$.  Since $\mu_i<d_i$ for all $i$, the monomials $LM(\varphi_jQ_i)=y_j^{d_j}y_i^{\mu_i}$ and $LM(\varphi_iQ_j)=y_i^{d_i}y_j^{\mu_j}$ are not equal in our monomial ordering.  Thus $LM(S(\varphi_i,\varphi_j)\ge \max(LM(\varphi_j Q_i), LM(\varphi_iQ_j))$.  It follows from Buchberger's Algorithm that $\{\varphi_1,\dots,\varphi_n\}$ is a Gr\"obner basis.
\qed

We will use the Multivariable Division Algorithm, by which any polynomial $g\in{\Bbb C}[y_1,\dots,y_n]$ may be written $g=A_1\varphi_1 + \cdots + A_n\varphi_n + R$  where $LM(g)\ge LM(A_j\varphi_j)$ for all $1\le j\le n$, and $R$ contains no terms divisible by any $LM(\varphi_j)$.  An important property of a Gr\"obner basis,  is that $g$ belongs to the ideal $\langle\varphi_1,\dots,\varphi_n\rangle$ if and only if $R=0$ (see, for instance, [CLO] or  [BW]).

If all fixed points have the same value of $\lambda$ as multiplier, then it follows that $\Phi$ must vanish on the whole zero set $Z(\varphi_1,\dots,\varphi_n)$.  Since we have a Gr\"obner basis, we easily determine the following:
\proclaim Corollary 4.4.  $\Phi\notin\langle \varphi_1,\dots,\varphi_n\rangle$.

\noindent{\it Proof. }  The leading monomial of $\Phi$ is $y_1^{d_1-1}\cdots y_n^{d_n-1}$, but this is not divisible by any of the leading monomials $LM(\varphi_j)=y_j^{d_j}$.  Since $\{\varphi_1,\dots,\varphi_n\}$ is a Gr\" obner basis, it follows that $\Phi$ does not belong to the ideal $\langle\varphi_1,\dots,\varphi_n\rangle$.
\qed

\medskip\noindent{\bf \S5.  Proof of the Theorem. }  In this section we prove the  Theorem, which will follow from Lemma 3.4, in combination with:  
\proclaim Proposition 5.1.  Suppose $F=f_n\circ\cdots\circ f_1$, $n\ge3$, is a composition of generalized H\'enon maps with $|\delta|<1$.  Suppose that $F$ has $d=d_1\cdots d_n$ distinct fixed points.  It is not possible that $d-1$ of these points have the same multipliers.

\noindent{\it Proof that Proposition 5.1 implies the Theorem. }  To prove the Theorem, it remains to deal with the case $|\delta|<1$.  If $f=f_1$ is a single generalized H\'enon map, we consider $F=f_1\circ f_1\circ f_1$ with $n=3$ and the same Julia set.  
Lemma 3.4 asserts that if  $J^+$ is $C^1$, there are $d-1$ saddle points with unstable multiplier $\lambda=d$.  So by Proposition 4.1 we conclude that $J^+$ cannot be $C^1$ smooth.  \qed

We give the proof of Proposition 5.1 at the end of this Section.  
For $J\subset\{1,\dots,n\}$, we write 
$$\Lambda_J:=\{(p')^L: L\subset J, |L|=|J|-2k, {\rm \ for\ some\ }, 1\le k\le |J| / 2\},$$
 We let $H_J$ denote the linear span of $\Lambda_J$.   To compare with our earlier notation, we note that $H_J\subset H_{\bf |J|-2}$ and that $(p')^J\notin H_J$.   The elements of $H_J$ depend only on the variables $y_j$ for $j\in J$.   
 Now we formulate a result for dividing certain terms by $\varphi_j$:
\proclaim Lemma 5.2.  Suppose that $J\subset\{1,\dots,n\}$ and $h\in H_J$.  Then for each $j\in J$ and $\alpha\in{\Bbb C}$, we have
$$(y_j-\alpha)\left( (p')^J + h\right) = A(y)\varphi_j + B(y) \left((p')^{J-\{j\}} + \rho_1 \right) +  (y_j-\alpha)\cdot \rho_2,\eqno(\dag)$$
where $\rho_1,\rho_2\in H_{J-\{j\}}$, and  $B=\eta_j(y_j) + d_j y_{j+1} + d_j\delta_j y_{j-1} $ with
$$\eta_j(y_j)  = y_jq_j'(y_j) - \alpha p'_j(y_j) - d_j q_j(y_j). \eqno(\ddag)$$ 
The leading monomials satisfy:
$$LM \left ((y_j-\alpha) \left( (p')^J+h \right) \right)=LM(A(y)\varphi_j)$$

%   $LT\left[  (y_j-\alpha)( (p')^J+h)\right] \ge LT(( (p')^{J-\{j\}} + K_3)\varphi_j)$
%
%\item{(1)}   $LT\left[  (y_j-\alpha)( (p')^J+h)\right] \ge LT(A\cdot \varphi_j)$
%\item{(2)}    $K_1 = (p')^{J-\{j\}} + K_3$
%\item{(3)}  $K_2,K_3\in H_{J-\{j\}}$
%\item{(4)} 
%\item{(5)}  $B = y_jq_j'(y_j) - \alpha p'_j(y_j) - q_j(y_j) + y_{j+1} + \delta_j y_{j-1}$
%\item{(6)}  $A=d_1K_1$

\noindent{\it Proof. }  Let us start with the case $J=\{1,\dots,m\}$, $m\le n$, and $j=1$, so $J-\{j\}=J_{\hat 1}=\{2,\dots,n\}$.  We divide by $p'_1$ and remove any factor of $p'_1$ in $h$.  This gives
$$(p')^J + h = p'_1(y_1)\mu_1 + \rho_2$$
where $\mu_1=(p')^{J_{\hat 1}} +\rho_1$, and $\rho_1,\rho_2\in H_{\{2,\dots,m\}}$, and  $\mu_1$, $\rho_1$, $\rho_2$ are independent of the variable $y_1$.  Thus
$$\eqalign{   (y_1-\alpha)\left (p')^J + h\right) & = (y_1-\alpha)(d_1 y_1^{d_1-1} + q'_1(y_1)) \mu_1 + (y_1-\alpha)\rho_2 \cr
& = d_1 y_1^{d_1}\mu_1 + (y_1 q'_1(y_1) - \alpha p_1'(y_1))\mu_1 + (y_1-\alpha)\rho_2\cr
& = (d_1\mu_1)\varphi_1 + (\eta_1(y_1) + d_1 y_2 + d_1\delta_1 y_n) \mu_1 + (y_1-\alpha) \rho_2}$$
where in the last line we substitute $\eta_1$ defined by $(\ddag)$.  Using $(*)$, we see that this gives $(\dag)$.  

It remains to look at the leading terms of $T_1:= (y_1-\alpha)\left ( (p')^J + h) \right)$ and $T_2:= d_1\mu_1\varphi_1$.  We see that $T_1$ and $T_2$ both contain nonzero multiples of $y_j \prod_{i=1}^m y_i^{d_i -1}$, and all other monomials in $T_1$ and $T_2$ have lower degree.  Thus we have $LM(T_1)=LM(T_2)$  for the graded ordering, independent of any ordering on the variables $y_1,\dots,y_n$.   The choices of $J=\{1,\dots,m\}$ and $j=1$ just correspond to a permutation of variables, and this does not affect the conclusion that $LM(T_1)=LM(T_2)$.  \qed

%As we will see in the proof of Proposition 4.1, it follows from Corollary 4.6 that $\lambda$ cannot be a multiplier for all fixed points of $f$.  If all of the fixed points have multiplier $\lambda$, with one exception $(\alpha,\beta)$, then $(y_1-\alpha)\Phi$ will vanish on $Z(\varphi_1,\dots,\varphi_n)$.  We consider this case next:
\proclaim Lemma 5.3.    For any $\alpha\in{\Bbb C}$, $(y_1-\alpha)\Phi\notin\langle \varphi_1,\dots,\varphi_n\rangle$.

\noindent{\it Proof. }   By [FM], we may assume that $p_j(y_j)=y_j^{d_j} + q_j(y_j) $, and ${\rm deg}(q_j)\le d_j-2$.  We consider two cases.  The first case is that there is at least one $j$ such that $\eta_j$ is not the zero polynomial.  If we conjugate by $f_{j-1}\circ\cdots\circ f_1$, we may ``rotate'' the maps in $f$ so that the factor $f_j$ becomes the first factor.  If there exists a $j$ for which $\eta_j(y_j)$ is non constant, we choose this for $f_1$.  Otherwise, if all the $\eta_j$ are constant, we choose $f_1$ to be any factor such that $\eta_1\ne0$.

We will apply the Multivariate Division Algorithm on $(y_1-\alpha)\Phi$ with respect to the set $\{\varphi_1,\dots,\varphi_n\}$.  We will find that there is a nonzero remainder, and since  $\{\varphi_1,\dots,\varphi_n\}$ is a Gr\"obner basis, it will follow that $(y_1-\alpha)\Phi$ does not belong to the ideal $\langle\varphi_1,\dots,\varphi_n\rangle$.   

We start with Lemma 4.2, according to which
$\Phi = p_1'\cdots p_n' + h$, where $h\in H_{\bf n-2} = H_{\{1,\dots,n\}}$.  The leading monomial of $(y_1-\alpha)\Phi$ is $y_1^{d_1}\prod_{i=2}^n y_i^{d_i-1}$, and $\varphi_1$ is the only element of the basis whose leading monomial divides this.  Thus we apply Lemma 5.2, with $J=\{1,\dots,n\}$, $j=1$, and $J_{\hat 1}:= J-\{j\}=\{2,\dots,n\}$.  This gives
$$\eqalign{ (y_1-\alpha)\Phi  &  = A_1 \varphi_1 + (\eta_1(y_1) + d_1 y_2 + d_1\delta_1 y_n) \left (\prod_{i=2}^n p_i'(y_i) +\rho_1 \right ) + (y_1-\alpha)\rho_2\cr
 = A_1 & \varphi_1  +\left[ d_1y_2 \left ((p')^{J_{\hat 1}}+\rho_1 \right ) \right]+  \left[d_1\delta_1 y_n\left ((p')^{J_{\hat 1}}+\rho_1 \right ) \right]+  \left[ \eta_1\left ((p')^{J_{\hat 1}}+\rho_1 \right )\right] + \ell.o.t \cr
  =A_1 & \varphi_1 + T_2 + T_n + R_1 + \ell.o.t}$$
where $\rho_1,\rho_2\in H_{\{2,\dots,n\}}$.  In particular, $T_2$ and $T_n$ depend on $y_2,\dots,y_n$ but not on $y_1$.  We note that $T_2$ (respectively,  $T_n$) contains a term divisible by $LM(\varphi_2)$ (respectively, $LM(\varphi_n)$).  We view  $R_1$ as a remainder term, and note that $LM(R_1)$ is divisible by $y_2^{d_2-1}\cdots y_n^{d_n-1}$, as well as the largest power of $y_1$ in $\eta_1(y_1)$.  By ``$\ell.o.t.$'' we mean that none of its monomials is divisible by $LM(R_1)$ or by any of the $LM(\varphi_j)$.

Now we apply Lemma 5.2 to $T_2$, this time with $J=\{2,\dots,n\}$ and $j=2$, with $J-\{2\} = J_{\hat 1\hat 2} = \{3,\dots,n\}$.  We have
$$\eqalign{ T_2 =  &A_2\varphi_2 + d_2y_3((p')^{J_{\hat 1\hat 2}}+ \rho_1^{(2)}) + d_2\delta_2 y_1 \left ( (p')^{J_{\hat 1\hat 2} } + \rho_1^{(2)} \right) + \eta_2(y_2) (p')^{J_{\hat 1\hat 2}} + \ell.o.t.\cr
= &A_2\varphi_2 + T_2^{(2)} + R_1^{(2)} + R_2^{(2)} + \ell.o.t. }$$
We see that $T_2^{(2)}$ contains terms that are divisible by $LM(\varphi_3)$, but the monomials in $R_1^{(2)}$ and $R_2^{(2)}$ are not divisible by $LM(\varphi_i)$ for any $i$.  The remainder term here is $R_1^{(2)}+R_2^{(2)}$, and we observe that this cannot cancel the largest term in $R_1$.  This is because $LM(R_1^{(2)})$ lacks a factor of $y_2$, and $LM(R_2^{(2)})$ is equal to $y_3^{d_3-1}\cdots y_n^{d_n-1}$ times the largest power of $y_2$ in $\eta_2(y_2)$, and by $(\ddag)$, this power is no bigger than $d_2-1$.  If $\eta_1$ is not constant, then we see that  $LM(R_1)>LM(R_2^{(2)})$.  If $\eta_1$ is constant, then $\eta_2$ must be constant, too, and again we have $LM(R_1)>LM(R_2^{(2)})$.  Thus, with our earlier notation, $R_1^{(2)} + R_2^{(2)} = \ell.o.t.$

We do a similar procedure with $T_n$, $T_2^{(2)}$,  etc., and again find that the remainder term does not contain a multiple of the leading monomial of  $R_1$.  We see that each time we do this process, the size of the exponent $L$ decreases in the term $(p')^L$.   When we have $L=\emptyset$, there are no terms that can be divided by any $LM(\varphi_j)$.  Thus we end up with
$$(y_1-\alpha)\Phi=A_1\varphi_1 + \cdots + A_n\varphi_n + R_1 + \ell.o.t.$$
and $LT((y_1-\alpha)\Phi)\ge LT(A_j\varphi_j)$ for all $1\le j\le n$, and none of the remaining terms is divisible by any of the leading monomials of $\varphi_j$.  Thus we have now finished the Multivariate Division Algorithm, and we have a nonzero remainder.  Thus $(y_1-\alpha)\Phi$ does not belong to the ideal of the $\varphi_j$'s.

Now we turn to the second case, in which $\eta_j=0$ for all $j$.  By [FM], we may assume that ${\rm deg}(q_j)\le d_j-2$. It follows that $\alpha=0$ and  $q_j=0$.  Thus $p_j = y_j^{d_j}$ for all $1\le j\le n$, so $p_j' = d_j y_j^{d_j-1}$, and $H_J$ consists of linear combinations of products $(p')^I= y_{i_1}^{d_{i_1}-1}\cdots y_{i_k}^{d_{i_k}-1}$ for $I = \{i_1,\dots,i_k\}\subset J$, for even $k\le |J|-2$.  We will go through the multivariate division algorithm again.  The principle is the same as before, but the details are different; in the first case we needed $n\ge2$, and now we will need $n\ge3$.   

Again, it is only $\varphi_1$ which has a leading monomial which can divide some terms in $(y_1-\alpha)\Phi$.  As before, we apply Lemma 5.2 with $J=\{1,\dots,n\}$, $j=1$, and $J-\{1\}=J_{\hat 1}=\{2,\dots,n\}$.  The polynomial in $(\ddag)$ becomes $B=d_jy_{j+1} + d_j\delta_j y_{j-1}$, and we have:
$$\eqalign{  y_1\Phi & = A_1\varphi_1 + d_1y_2 \left( (p')^{J_{\hat 1}} + \rho_1\right) + d_1\delta_1 y_n\left( (p')^{J_{\hat 1}}+ \rho_1\right)  + y_1\rho_2\cr
& = A_1\varphi_1 + T_2 + T_n + \ell.o.t. }$$
where $\rho_1,\rho_2\in H_{\{2,\dots,n\}}$.  Now we apply Lemma 5.2 to divide $T_2$ (respectively $T_n$) by $\varphi_2$ (respectively $\varphi_n$).  This yields:
$$y_1\Phi=A_1\varphi_1+ A_2\varphi_2 + A_n\varphi_n + T_3+T_n+  R+\ell.o.t.$$
where
$$T_3= d_1d_2y_3\left( (p')^{J_{\hat 1\hat 2} }+ \tilde \rho_3 \right), \ \ \  T_n= d_1d_n\delta_1\delta_n y_{n-1} \left ( (p')^{J_{\hat 1\hat n}} + \tilde \rho_n \right)  $$
with $\tilde\rho_3\in H_{\{3,\dots,n\}}$ and $\tilde\rho_n\in H_{\{2,\dots,n-1\}}$, and 
$$R =  \left(d_1d_2\delta_2 y_1 y_n^{d_n-1}  + d_1d_n\delta_1 y_1 y_{2}^{d_2-1}  \right) \prod_{i=3}^{n-1} y_i^{d_i-1}$$
Since $n>2$,  $R$ is not the zero polynomial.  We will continue the Multivariate Division Algorithm by dividing $T_3$ by $\varphi_3$ and $T_n$ by $\varphi_n$, but we see that any terms created cannot cancel $R$.  Thus when we finish the Division Algorithm, we will have a nonzero remainder.  As in the previous case, we conclude that $y_1\Phi$ is not in the ideal $\langle\varphi_1,\dots,\varphi_n\rangle$.  
\qed

\noindent {\it Proof of Proposition 5.1.  }  The fixed points of $f$ coincide with the elements of $Z(\varphi_1,\dots,\varphi_n)$, which is a variety of pure dimension zero.  Saddle points have multiplicity 1, and since there are $d-1$ of these, and since the total multiplicity is $d$,  there must be one more fixed point,  also of multiplicity 1.    It follows that the ideal $I:=\langle \varphi_1,\dots,\varphi_n\rangle$ is equal to its radical (see [BW]).  Since the saddle points all have multiplier $\lambda$, $\Phi$ must vanish at all the saddle points.  If $(\alpha,\beta)$ is the other fixed point, we conclude that $(y_1-\alpha)\Phi$ vanishes at all the fixed points.  Thus $(y_1-\alpha)\Phi$ belongs to the radical of $I$, and thus $I$ itself.  This contradicts Lemma 5.3, which completes the proof of Proposition 5.1.  
\qed

\medskip
\centerline{\bf Appendix:  Non-smoothness of $J$, $J^*$, and $K$}
\medskip
Let us turn our attention to other dynamical sets for polynomial diffeomorphisms of positive entropy.  These are $J:=J^+\cap J^-$, $K:= K^+\cap K^-$, and the  set $J^*$, which coincides with the closure of the set of periodic points of saddle type.  (See [BS1], [BS3], and [BLS] for other characterizations of $J^*$.) We have $J^*\subset J\subset K$.  We note that none of these sets can be a smooth 3-manifold: otherwise, for any saddle point $p$, it would be a bounded set containing $W^s(p)$ or $W^u(p)$, which is the holomorphic image of ${\Bbb C}$.  The following was suggested by Remark 5.9 of Cantat in [C]; we sketch his proof:  

\proclaim Proposition A.1.  If $J=J^*$, then it is not a smooth 2-manifold. 

\noindent{\it Proof.}    Let  $p$  be a saddle point, and let  $W^u(p)$ be the unstable manifold.  The slice $J \cap W^u(p)$  is smooth and invariant under multiplication by the multiplier of  $Df$.  This means that in fact, the multiplier must be real, and the  restriction of  $G^+$ to the slice must be linear on each (half-space) component of  $W^u(p)- J$. 

The identity  $G^+\circ f = d\cdot G^+$ means that the canonical metric (defined in [BS8]) is multiplied by  $d$.  Thus $f$ is quasi-expanding on  $J^*$.  Now, applying this argument to  $f^{-1}$ we get that $f$ is quasi-hyperbolic.  Further, $J^*=J$, so it is quasi-hyperbolic on $J$.  If  $f$ fails to be hyperbolic, then by  [BSm] there will be a one-sided saddle point, which can not happen since  $J$ is smooth.

Now that $f$ is hyperbolic on $J$, there is a splitting $E^s \oplus E^u$  of the tangent bundle, so we conclude that $J$ is a 2-torus.  The dynamical degree must be the spectral radius of an invertible 2-by-2 integer matrix, but this means it is not an integer, which contradicts the fact the the dynamical degree of a H\'enon map is its algebraic degree.
\qed

\proclaim Proposition A.2.  Suppose that the complex jacobian is not equal to $\pm1$.  Then for each saddle (periodic) point $p$ and each neighborhood $U$ of $p$, neither $J\cap U$ nor $J^*\cap U$ nor $K\cap U$ is a $C^1$ smooth 2-manifold.

\noindent{\it Proof.}   Let us write $M:=J\cap U$ and $g:=f|_M$.  (The following argument works, too, if we take $M=J^*\cap U$ or $M=K\cap U$.)   The tangent space $T_pM$ is invariant under $Df$.  The stable/unstable spaces $E^{s/u}\subset T_p{\Bbb C}^2$ are invariant under $D_pf$.  The space $E^s$ (or $E^u$) cannot coincide with $T_pM$, for otherwise the complex stable manifold $W^s(p)$ (or $W^u(p)$) would be locally contained in $M$, and thus globally contained in $J$.  But the $W^{s/u}$ are uniformized by ${\Bbb C}$, whereas $J$ is bounded.  We conclude that $p$ is a saddle point for $g$, and thus the local stable manifold $W^s_{\rm loc}(p; g)$ is a $C^1$-curve inside the complex stable manifold $W^s(p)$.  As in Lemma 3.3, we conclude that the multiplier for $D_pf|_{E^u_p}$ is $\pm d$ and the multiplier for $D_pf|_{E^s_p}$ is $\pm 1/d$.  Thus the complex Jacobian is $\delta=\pm 1$. 
\qed

\noindent{\bf Solenoids. } The two results above concern smoothness, but no example is known where $J$,  $J^*$ or $K$ is even a topological 2-manifold.  In the cases where $J^+$ has been shown to be a topological 3-manifold (see [FS], [HO2], [Bo] and [RT]) it also happens that $J$ is a (topological) real solenoid, and in these cases it is also the case that $J=J^*$.  Further, for every saddle (periodic) point $p$, there is a real arc $\gamma_p=W^u_{\rm loc}(p)\cap J$.  If we apply the argument of Proposition A.2 to this case, we conclude that $\gamma_p$ is not $C^1$ smooth.

\vfill\eject
\bigskip
\centerline{\bf References}
\medskip

\item{[BW]}  T. Becker and V. Weispfenning,  Gr\"obner Bases, A computational Approach to Commutative Algebra,  Springer-Verlag, Berlin and New York, 1993, xxii + 574 pp.

\item{[BLS]}  E. Bedford, M. Lyubich, and J. Smillie,  Polynomial diffeomorphisms of ${\bf C}^2$. IV. The measure of maximal entropy and laminar currents. Invent. Math. 112 (1993), no. 1, 77--125.

\item{[BS1]}  E. Bedford and J. Smillie, Polynomial diffeomorphisms of ${\bf C}^2$: currents, equilibrium measure and hyperbolicity. Invent. Math. 103 (1991), no. 1, 69--99. 

\item{[BS2]} E. Bedford and J. Smillie, Polynomial diffeomorphisms of ${\bf C}^2$. II. Stable manifolds and recurrence. J. Amer. Math. Soc. 4 (1991), no. 4, 657--679. 

\item{[BS3]}  E.  Bedford and J. Smillie,  Polynomial diffeomorphisms of ${\bf C}^2$. III. Ergodicity, exponents and entropy of the equilibrium measure. Math. Ann. 294 (1992), no. 3, 395--420.

\item{[BS8]}  E. Bedford and J. Smillie,  Polynomial diffeomorphisms of ${\bf C}^2$. VIII. Quasi-expansion. Amer. J. Math. 124 (2002), no. 2, 221--271.

\item{[BSm]} E. Bedford and J. Smillie,  Real polynomial diffeomorphisms with maximal entropy: Tangencies. Ann. of Math. (2) 160 (2004), no. 1, 1--26. 

\item{[Bo]} S. Bonnot,  Topological model for a class of complex H\'enon mappings, Comment. Math. Helv. 81 (2006), no. 4, 827--857.

\item{[C]}  S. Cantat,  Bers and H\'enon, Painlev\'e and Schr\"odinger.  Duke Math. Journal, vol 149 (2009), no. 3, 411--460 

\item{[CLO]}  D. Cox, J. Little, D. O'Shea,  Ideals, varieties, and algorithms.  Springer, New York, 2007. xvi+551 pp.

\item{[FS]} J.E. Forn\ae ss and N. Sibony,  Complex H\'enon mappings in ${\bf C}^2$ and Fatou-Bieberbach domains. Duke Math. J. 65 (1992), no. 2, 345--380. 

\item{[FM]}  S. Friedland and J. Milnor,  Dynamical properties of plane polynomial automorphisms.
Ergodic Theory Dynam. Systems 9 (1989), no. 1, 67--99. 

\item{[Ha]}  M. Hakim,  Attracting domains for semi-attractive transformations of ${\Bbb C}^p$. Publ. Mat. 38 (1994), no. 2, 479--499.

\item{[H]}  J.H. Hubbard,  The H\'enon mapping in the complex domain. Chaotic dynamics and fractals (Atlanta, Ga., 1985), 101--111, Notes Rep. Math. Sci. Engrg., 2, Academic Press, Orlando, FL, 1986.

\item{[HO1]}  J.H. Hubbard and R. Oberste-Vorth,  H\'enon mappings in the complex domain. I. The global topology of dynamical space. Inst. Hautes \'Etudes Sci. Publ. Math. No. 79 (1994), 5--46. 

\item{[HO2]}  J.H. Hubbard and R. Oberste-Vorth,  H\'enon mappings in the complex domain. II. Projective and inductive limits of polynomials. Real and complex dynamical systems (Hiller\o d, 1993), 89--132, NATO Adv. Sci. Inst. Ser. C Math. Phys. Sci., 464, Kluwer Acad. Publ., Dordrecht, 1995. 

\item{[M1]}  J. Milnor, Dynamics in one complex variable. Third edition. Annals of Mathematics Studies, 160. Princeton University Press, Princeton, NJ, 2006.

\item{[M2]}  J. Milnor,  Topology from the differentiable viewpoint. Based on notes by David W. Weaver. The University Press of Virginia, Charlottesville, Va. 1965 ix+65 pp. 

\item{[RT]}  R. Radu and R. Tanase,  A structure theorem for semi-parabolic H\'enon maps,  

arXiv:1411.3824

\item{[Sa]}  H.  Samelson,  Orientability of hypersurfaces in ${\Bbb R}^n$. Proc.\ Amer.\ Math.\ Soc.\ 22 (1969) 301--302. 

\item{[Sm]}  J. Smillie,  The entropy of polynomial diffeomorphisms of ${\Bbb C}^2$. Ergodic Theory Dynam. Systems 10 (1990), no. 4, 823--827. 

\item{[U]} T. Ueda,  Local structure of analytic transformations of two complex variables.~I.
{\sl J. Math. Kyoto Univ.} 26 (1986), no. 2, 233--261.

\bigskip
\rightline{Eric Bedford}

\rightline{Stony Brook University}

\rightline{Stony Brook, NY 11794}

\rightline{\tt ebedford@math.sunysb.edu}

\bigskip
\rightline{Kyounghee Kim}

\rightline{Florida State University}

\rightline{Tallahassee, FL 32306}

\rightline{\tt kim@math.fsu.edu}

\bye